\documentclass[oneside]{amsart}
\usepackage{amsmath,graphicx,amssymb}

\newtheorem{mainthm}{Theorem}
\newtheorem{thm}{Theorem}[section]
\newtheorem{lem}[thm]{Lemma}
\newtheorem{cor}[thm]{Corollary}

\theoremstyle{definition}
\newtheorem{exmp}[thm]{Example}
\newtheorem{defn}[thm]{Definition}
\newtheorem{rem}[thm]{Remark}

\newcommand{\blackboard}[1]{\ensuremath{\mathbb{#1}}}
\newcommand{\script}[1]{\ensuremath{\mathcal{#1}}}
\newcommand{\smallcaps}[1]{\textrm{\textsc{#1}}}

\newcommand{\N}{\blackboard{N}}
\newcommand{\Z}{\blackboard{Z}}

\newcommand{\R}{\blackboard{R}}

\newcommand{\euler}{\raisebox{.5mm}{\ensuremath{\chi}}}

\newcommand{\psig}{{\mbox{\rm P}\Sigma}}

\newcommand{\vcd}{{\mbox{\rm vcd}}}
\newcommand{\val}{{\mbox{\rm val}}}
\newcommand{\Aut}{\smallcaps{Aut}}
\newcommand{\Out}{\smallcaps{Out}}
\newcommand{\Inn}{\smallcaps{Inn}}
\newcommand{\OWh}{\smallcaps{OWh}}
\newcommand{\Wh}{\smallcaps{Wh}}
\newcommand{\OFR}{\smallcaps{OFR}}
\newcommand{\FR}{\smallcaps{FR}}

\newcommand{\MM}{{\mbox{MM}}}

\newcommand{\es}{\smallcaps{EdgeSizes}}
\newcommand{\blocks}{\smallcaps{Blocks}}
\newcommand{\vs}{\smallcaps{Vertices}}
\newcommand{\edges}{{\smallcaps{Edges}}}
\newcommand{\hyper}{\smallcaps{HT}}

\newcommand{\PHF}{\smallcaps{PHF}}
\newcommand{\RHT}{\smallcaps{RHT}}
\newcommand{\wt}{{\smallcaps{Weight}}}
\newcommand{\ewt}{{\smallcaps{E-Weight}}}
\newcommand{\vwt}{{\smallcaps{V-Weight}}}
\newcommand{\rk}{{\mbox{\rm rank}}}

\begin{document}

\title[$\euler(\Wh(G_1 \ast \cdots \ast G_n)) 
= \euler(G_1 \ast \cdots \ast G_n)^{n-1}$] 
{The Euler characteristic of  the Whitehead\\ automorphism group of a free product}

\author[C.~Jensen]{Craig Jensen$\ \!{ }^1$} 
\address{Dept. of Mathematics\\
         University of New Orleans\\
         New Orleans, LA 70148}
\email{jensen@math.uno.edu}

\author[J.~McCammond]{Jon McCammond$\ \!{ }^2$} 
\address{Dept. of Mathematics\\ 
         University of California\\ 
         Santa Barbara, CA 93106}
\email{jon.mccammond@math.ucsb.edu}

\author[J.~Meier]{John Meier$\ \!{ }^3$}
\address{Dept. of Mathematics\\
         Lafayette College\\
         Easton, PA 18042}
\email{meierj@lafayette.edu}


\begin{abstract}
A combinatorial summation identity over the lattice of labelled
hypertrees is established that allows one to gain concrete information
on the Euler characteristics of various automorphism groups of free
products of groups.  In particular, we establish formulae for the
Euler characteristics of: the group of Whitehead automorphisms
$\Wh(\ast_{i=1}^n G_i)$ when the $G_i$ are of finite homological type;
$\Aut(\ast_{i=1}^n G_i)$ and $\Out(\ast_{i=1}^n G_i)$ when the $G_i$
are finite; and the palindromic automorphism groups of finite rank
free groups.
\end{abstract}

\footnotetext[1]{Partially supported by Louisiana Board of Regents RCS
  contract no.\ LEQSF-RD-A-39}
\footnotetext[2]{Partially supported by NSF grant no.\ DMS-0101506}
\footnotetext[3]{Partially supported by an AMS Centennial Research Fellowship}
\maketitle

\section{Introduction}\label{intro}
Let $G = G_1\ast \cdots \ast G_n$, where the $G_i$ are non-trivial groups.
There are various subgroups of $\Aut(G)$ and $\Out(G)$ (such as the
Whitehead automorphism group) that are only defined with respect to a
specified free product decomposition of $G$.  In this article we
calculate the Euler characteristics of several such groups.
Postponing definitions for the moment, our main result is the
following.

\begin{mainthm}\label{mainthm:owh}
If $G = G_1 \ast\cdots \ast G_n$ is a free product of groups where
$\euler(G)$ is defined, then the Euler characteristic of the group of
outer Whitehead automorphisms is
\[
\euler(\OWh(G))=\euler(G)^{n-2} = \left[ \euler(G_1) + \cdots +
  \euler(G_n) - (n-1)\right]^{n-2}\ .
\]
As a consequence, the Euler characteristic of the group of Whitehead
automorphisms is $\euler(\Wh(G))=\euler(G)^{n-1}$.
\end{mainthm}

A natural situation to consider is when each group $G_i$ is finite.
In this case the $\vcd$ of $\Aut(G)$ is $n-1$ (independently
established in \cite{KrVo93} and \cite{McMi96}) and $\Wh(G)$ is a
finite index subgroup of $\Aut(G)$.  Theorem~\ref{mainthm:owh} thus
implies:

\begin{mainthm}\label{mainthm:aut}
If $G = G_1 \ast \cdots \ast G_n$ is a free product of finite groups and
$\Omega \subset \Sigma_n$ is the subgroup of the symmetric group given
by permuting isomorphic groups $G_i$, then the Euler characteristics
of the Fouxe-Rabinovitch automorphism group of $G$, the full
automorphism group of $G$ and the outer automorphism group of $G$ are:
\[\begin{array}{rcl}
\euler(\FR(G)) &=& \euler(G)^{n-1}  \displaystyle\prod_{i=1}^n |\Inn(G_i)|\\
\euler(\Aut(G)) &=& \euler(G)^{n-1}
        |\Omega|^{-1} \displaystyle\prod_{i=1}^n |\Out(G_i)|^{-1}\\
\euler(\Out(G)) &=& \euler(G)^{n-2} 
        |\Omega|^{-1} \displaystyle\prod_{i=1}^n |\Out(G_i)|^{-1}\\
\end{array}\]
\end{mainthm}

The automorphism groups mentioned in these theorems are defined as
follows.

\begin{defn}[Groups of automorphisms]\label{def:aut-gps}
Let $G_1 \ast \cdots \ast G_n$ be a free product of groups.  As usual, we
use $\Aut(G)$ to denote the \emph{full automorphism group} of $G$ and
$\Out(G)$ for the group of \emph{outer automorphisms}.  If 
$g_i \in G_i \setminus \{1\}$ let $\alpha_j^{g_i}$ denote
the automorphism described by
\[
\alpha_j^{g_i}(g) =  \left\{
\begin{array}{cc}
g & g \in G_k\mbox{ where }k\ne j\cr
g^{g_i} & g \in G_k\mbox{ where }k = j\cr
\end{array}  \right.\ .
\]
The \emph{Fouxe-Rabinovitch group} $\FR(G)$ is the subgroup of
$\Aut(G)$ generated by $\{\alpha_j^{g_i}~|~i \ne j\}$.  The group of
\emph{Whitehead automorphisms} is generated by all the
$\alpha_j^{g_i}$ (including the possibility that $i=j$). 
In general, the \emph{symmetric automorphisms} of a free
product consist of those automorphisms that send each
factor group to a conjugate of a (possibly different) factor. It is
well-known that the group $\Wh(G)$ is the kernel of the map
$f:\Sigma\Aut(G)\to \Out(\prod_i G_i)$ and $\FR(G)$ is the kernel of the map
$g:\Sigma\Aut(G) \to \Aut(\prod_i G_i)$, where $\Sigma\Aut(G)$ denotes 
the group of symmetric automorphisms.

A quick chasing of definitions shows that the images of $\FR(G)$ and
$\Wh(G)$ in $\Out(G)$ are the same --- $\OFR(G)=\OWh(G)$ --- and so
$\OFR(G)$ will not appear again.  (For more background on these
definitions see \cite{McMi96}, especially page 48.)
\end{defn}

A {\em palindromic automorphism} of a free group
$F_n=F(a_1,\ldots,a_n)$ is an automorphism sending each generator
$a_i$ to a palindromic word.  The {\em elementary palindromic
automorphisms} are those which send each $a_i$ to a palindrome of
odd length 
with $a_i$ as its centermost letter.  Lastly, the {\em pure
palindromic automorphisms} are those which send each $a_i$ to an odd
palindrome with either $a_i$ or $a_i^{-1}$ in the center.  After
identifying these automorphisms of $F_n$ with ones of $\ast_{i=0}^n
\Z_2$, we are able to prove:

\begin{mainthm}\label{mainthm:palindromic}
If $F_n$ is a free group of rank $n$, then the Euler characteristics
of the elementary palindromic automorphism group of $F_n$, the pure
palindromic automorphism group of $F_n$, and the palindromic
automorphism group of $F_n$ are:
\begin{eqnarray*}
\euler(E\Pi A(F_n)) & = & (1-n)^{n-1}\\
\euler(P\Pi A(F_n)) & = & \frac{(1-n)^{n-1}}{2^n} \\
\euler(\Pi A(F_n))  & = & \frac{(1-n)^{n-1}}{2^n \cdot n!}
\end{eqnarray*}
\end{mainthm}

\vspace{10pt}

Even though our stated theorems are related to group cohomology, the
underlying arguments are primarily combinatorial.  In the next section
we give background information on polynomial identities and a
particularly useful identity relating a multi-variable polynomial with
the number of rooted trees and planted forests.
Section~\ref{sec:partitions} is the heart of the combinatorial work,
where we use the polynomial identities to establish a number of
partition identities and ultimately an identity involving a summation
of particular weights over the set of hypertrees
(Theorem~\ref{thm:comb-main}).  Then in Section~\ref{sec:euler}, the
notion of an Euler characteristic of a group is reviewed, the space
constructed by McCullough and Miller is outlined, and an approach to
calculating the Euler characteristic of $\OWh(G)$ is sketched.  In
Section~\ref{sec:main}, the hypertree identity is used to complete the
calculation and we complete the proofs of Theorems~\ref{mainthm:owh}
and \ref{mainthm:aut}, as well as give some concrete examples.  In
Section~\ref{sec:palindromic} we turn to the group of palindromic
automorphisms of $F_n$ and establish
Theorem~\ref{mainthm:palindromic}.

Although the current ordering of the sections in the paper presents
results in the order in which they are proven, an alternate ---
possibly more conceptual and motivated --- ordering would be
\[
\mbox{\S 1} \Rightarrow \mbox{\S 4} 
\Rightarrow \mbox{\S 5} \Rightarrow \mbox{\S 6}
\Rightarrow \mbox{\S 3} \Rightarrow \mbox{\S 2} \ .
\]  
Sections 4--6 show how the Euler characteristic calculations can be
derived from a particular weighted sum over hypertrees.  This
hypertree identity is reduced to two weighted sums over the set of
partitions in \S 3, and finally, the necessary partition identities
are established using a polynomial identity described in \S 2.

\section{Polynomials and Trees}\label{sec:poly}
In this section we record some well-known polynomial identities that
we need later in the article.  The first two identities were
discovered by Niels Abel. In order to make them easier to state we
introduce the following function.

\begin{defn}[An interesting basis]\label{def:basis}
For each non-negative integer $k$ let $h_k(x)$ denote the polynomial
$x(x+k)^{k-1}$.  Note that when $k=0$ the expression $x(x+0)^{-1}$ is
not quite the same as the constant polynomial $1$ since it is
undefined at $x=0$, but this situation can be remedied by explicitly
defining $h_0(x) =1$.  The polynomials $\{h_k(x)\}_{k\geq 0}$ form a
basis for the ring of polynomials in $x$, because the degree of
$h_k(x)$ is $k$ and thus there is exactly one polynomial of each
degree.  As a consequence, the ``monomials'' of the form
$h_{k_1}(x_1)\cdot h_{k_2}(x_2) \cdots h_{k_n}(x_n)$ with each
$k_i\geq 0$ form a basis for the ring of polynomials in commuting
variables $x_1,x_2,\ldots, x_n$.
\end{defn}

\begin{lem}[Abel's identities]\label{lem:abel}
For each non-negative integer $n$, the two polynomial identities
listed below are valid.
\begin{equation}\label{eq:abel1}
(x+y+n)^n = \sum_{i+j=n} \binom{n}{i}\ h_i(x)\cdot (y+j)^{j}
\end{equation}
\begin{equation}\label{eq:abel2}
h_n(x+y) = \sum_{i+j=n} \binom{n}{i}\ h_i(x)\cdot h_j(y)
\end{equation}
\end{lem}

These identities are readily established by induction (see for example
\cite{Ro84}).  The second identity listed has an obvious
multi-variable generalization using multinomial coefficients.

\begin{rem}
The exponential generating function for $h_n(x)$ is the function
\cite{Concrete} calls $(\script{E}_1(z))^x$.  The generating function
equivalent of Lemma~\ref{lem:abel} is the fact that
$(\script{E}_1(z))^x \cdot (\script{E}_1(z))^y =
(\script{E}_1(z))^{x+y}$.
\end{rem}

\begin{figure}[hptb]
\includegraphics[width=2.0in]{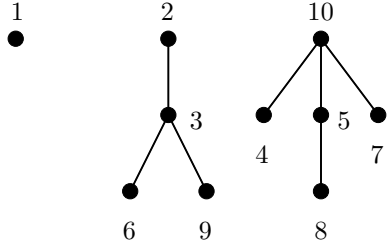}
\caption{This planted forest on $[10]$ has roots $1$, $2$ and $10$.
Its degree sequence is $(0,1,2,0,1,0,0,0,0,3)$ and the associated
degree sequence monomial is $x_2 x_3^2 x_5 x_{10}^3$.}
 \label{fig:plantedforest}
 \end{figure}

\begin{defn}[Trees and forests]
A \emph{forest} is a simplicial graph with no cycles and a \emph{tree}
is a connected forest.  All of the trees and forests considered here
have distinguishable vertices.  In particular, if $[n]:=\{1,\ldots,
n\}$ then it is common practice to consider the trees or forests that
can be described using $[n]$ as the vertex set.  A tree or forest with
a distinguished vertex selected in each connected component is said to
be \emph{rooted} or \emph{planted}, respectively.  In a planted forest
we consider each edge as oriented away from the root in its component.
The \emph{rooted degree} of a vertex $i$ in a planted forest is the
number of edges starting at $i$.  Notice that the rooted degree of a vertex
is one less than its valence unless that vertex happens to be a root,
in which case the two numbers are equal.  The \emph{degree sequence}
of a planted forest on $[n]$ is the $n$-tuple of rooted degrees
$\delta_i$ (in order) for $i \in [n]$, and the associated \emph{degree
sequence monomial} is $x_1^{\delta_1} x_2^{\delta_2} \cdots
x_n^{\delta_n}$.  See Figure~\ref{fig:plantedforest} for an example.
\end{defn}

\begin{thm}[A forest identity, Theorem 5.3.4 in \cite{StV2}]\label{thm:forest}
Let $\delta = (\delta_1,\ldots,\delta_n) \in \N^n$ with $\sum \delta_i
= n-k$.  The number $N(\delta)$ of planted forests $\sigma$ on the
vertex set $[n]$ (necessarily with $k$ components) with ordered degree
sequence $\Delta(\sigma) = \delta$ is given by \[ N(\delta) =
\binom{n-1}{k-1} \binom{n-k}{\delta_1, \delta_2, \ldots, \delta_n}.\]
Equivalently, \[ \sum_\sigma x_1^{\deg 1} \cdots x_n^{\deg n} =
\binom{n-1}{k-1} (x_1+ \cdots + x_n)^{n-k},\] where $\sigma$ ranges
over all planted forests on $[n]$ with $k$ components.
\end{thm}

As an illustration, consider the $6$ planted forests on $[3] :=
\{1,2,3\}$ with exactly $2$ components.  Each such forest contains a
single edge and a unique vertex with a nonzero degree so that the
final sum in the theorem is $2(x_1 + x_2 + x_3)$ as claimed.

\section{Partitions and Hypertrees}\label{sec:partitions}

In this section we establish several identities involving summations
over sets of partitions and hypertrees.  Our Euler characteristic
computations ultimately involve sums over the poset of labelled
hypertrees, a poset whose internal structure is quite similar to the
well-studied poset of partitions.  After certain polynomial identities
over the poset of partitions are established in the first subsection,
we convert them into hypertree identities in the second subsection.

\subsection{Partition Identities}

\begin{defn}[Partitions]
Recall that a \emph{partition} $\sigma$ of $[n]$ is a collection of
pairwise disjoint subsets of $[n]$ whose union is $[n]$.  These
subsets (i.e. the elements of $\sigma$) are called \emph{blocks}.  The
collection of all partitions of $[n]$ is denoted $\Pi_n$ and the
subcollection of all partitions with exactly $k$ blocks is denoted
$\Pi_n^k$.
\end{defn}

\begin{defn}[Partition sums]
For each $i\in \{1,\ldots, n\}$, let $w_i$ be a fixed choice of a
symmetric polynomial in $i$ variables.  These functions can be used to
define weights on the blocks of a partition as follows.  Given a
partition $\sigma$ of $[n]$ and an element $\tau$ in $\sigma$ of size
$i$ we evaluate $w_i$ at the variables $x_j$, $j\in \tau$.  The
corresponding \emph{partition weight function} is given by taking the
product over all the blocks in the partition $\sigma$.  The
\emph{partition lattice function} $f(n,k)$ is the sum of the partition
weight functions over all partitions of $[n]$ with exactly $k$ blocks.
\end{defn}

One particular symmetric polynomial occurs often enough below so as to
merit a special notation.  For any set $\tau \subset [n]$ define
$s_\tau$ as the sum of the variables $x_i$ over all $i\in \tau$.  For
example, $s_{\{1,2,4\}} = x_1+x_2+x_4$.

\begin{exmp}
Let $n=3$ and $k=2$.  If the symmetric polynomials $s_\tau$ are used
as block weights, then the partition lattice function $f(3,2)$ is
$x_1(x_2+x_3) + x_2(x_1+x_3) + x_3(x_1+x_2) = 2(x_1x_2 + x_1x_3 +
x_2x_3)$.
\end{exmp}

Notice that it is clear from the definition that $f(n,k)$ is itself
always a symmetric function on $n$ variables, but it is not at all
clear whether one should expect closed form formulae for $f(n,k)$ for
a particular choice of block weights.  Using the forest identity from
the previous section, it is relatively easy to establish at least one
such summation.

\begin{lem}[$1^{st}$ partition identity]\label{lem:part1}
The partition weight function $f(n,k)$ that weights each block $\tau$
by $(s_\tau)^{|\tau|-1}$ has closed form $\binom{n-1}{k-1}
(s_{[n]})^{n-k}$.  In other words,
\[f(n,k) = \sum_{\sigma \in \Pi_n^k} \left( \prod_{\tau\in
  \blocks(\sigma)} (s_\tau)^{|\tau|-1}\right) = \binom{n-1}{k-1}
    (s_{[n]})^{n-k}.\]
\end{lem}

\begin{proof}
Using the forest identity (Theorem~\ref{thm:forest}), the weight on
block $\tau$ is equal to a sum of  degree sequence
monomials for  all
rooted trees on the underlying set $\tau$.  Replacing each block
weight with this summation over rooted trees reveals that the left
hand side is actually a sum of the degree sequence monomials over all
planted forests with exactly $k$ components.  A second use of
Theorem~\ref{thm:forest} completes the proof.
\end{proof}

By specializing each variable $x_i$ to $1$, we obtain the following
result as an immediate corollary.

\begin{cor}[$2^{nd}$ partition identity]\label{cor:part2}
The partition weight function $f(n,k)$ which weights each block $\tau$
by $|\tau|^{|\tau|-1}$ has closed form $\binom{n-1}{k-1}\ n^{n-k}$.
In other words,
\[f(n,k) = \sum_{\sigma \in \Pi_n^k} \left( \prod_{\tau\in
  \blocks(\sigma)} |\tau|^{|\tau|-1}\right) = \binom{n-1}{k-1}\
    n^{n-k}.\]
\end{cor}

A second corollary uses the alternative basis described in
Definition~\ref{def:basis}.  Recall that $h_k(x)$ denotes the
polynomial $x(x+k)^{k-1}$ so that if $\tau = \{1,2,4,5\}$ then
$h_{|\tau|-1}(s_\tau) = (x_1+x_2+x_4+x_5) (x_1+x_2+x_4+x_5 +3)^2$.

\begin{cor}[$3^{rd}$ partition identity]\label{cor:part3}
The partition weight function $f(n,k)$ which weights each block $\tau$
by $h_{|\tau|-1}(s_\tau)$ has closed form $\binom{n-1}{k-1}
h_{n-k}(s_{[n]})$.  In other words,
\[f(n,k) = \sum_{\sigma \in \Pi_n^k} \left( \prod_{\tau\in
  \blocks(\sigma)} h_{|\tau|-1}(s_\tau) \right)= \binom{n-1}{k-1}
  h_{n-k}(s_{[n]}).\]
\end{cor}

\begin{proof}
Consider the linear operator $T$ on $\R[x_1,\ldots,x_n]$ (viewed as an
$\R$-vector space) that sends the monomial $x_1^{\delta_1}\cdots
x_n^{\delta_n}$ to the ``monomial'' $h_{\delta_1}(x_1)\cdots
h_{\delta_n}(x_n)$.  If both sides of the equation in the statement of
Lemma~\ref{lem:part1} are expanded using the usual binomial theorem,
then the fact that the coefficients of the corresponding monomials on
each side are equal is essentially the content of the first part of
the forest identity from Theorem~\ref{thm:forest}.  If each monomial on
each side is now replaced with the corresponding $h$-monomial
(i.e. its image under $T$) --- without changing the coefficients used
--- then the $h$-version of the binomial theorem (Lemma~\ref{lem:abel})
can be used to recombine each side so as to establish the equality
claimed above.
\end{proof}

\subsection{Hypertree Identities}

\begin{defn}[Hypertrees and hyperforests]\label{def:hypertrees}
A \emph{hypergraph} $\Gamma$ is an ordered pair $(V,E)$ where $V$ is
the set of vertices and $E$ is a collection of subsets of $V$---called 
\emph{hyperedges}---each containing at least two elements.  A
\emph{walk} in a hypergraph $\Gamma$ is a sequence
$v_0,e_1,v_1,\ldots,v_{n-1},e_n,v_n$ where for all $i$, $v_i \in V$,
$e_i \in E$ and for each $e_i$, $\{v_{i-1},v_i\}\subset e_i$.  A
hypergraph is \emph{connected} if every pair of vertices is joined by a
walk.  A \emph{simple cycle} is a walk that contains at least two
edges, all of the $e_i$ are distinct and all of the $v_i$ are distinct
except $v_0=v_n$.  A hypergraph with no simple cycles is a
\emph{hyperforest} and a connected hyperforest is a \emph{hypertree}.
Note that the no simple cycle condition implies that distinct edges in
$\Gamma$ have at most one vertex in common.  All of the hypertrees and
hyperforests considered here have distinguishable vertices.  In
particular, it is common practice to mainly consider the hypertrees or
hyperforests that can be described using $[n]$ as the vertex set.
Examples of hypertrees on $[n]$ are shown in
Figures~\ref{fig:hypercorresp} and \ref{fig:interval}, where
hyperedges are indicated by polygons.
\end{defn}

Many of the concepts previously defined for trees and forests have
natural extensions to hypertrees and hyperforests.  For example, a
hypertree or hyperforest with a distinguished vertex selected in each
connected component is said to be \emph{rooted} or \emph{planted},
respectively.  In a planted hyperforest it is common practice to
consider each hyperedge to be ``oriented'' by gravity,  away from the root of the
component.  In this case the ``orientation'' of a hyperedge
distinguishes between the unique vertex closest to the root and all of
the remaining vertices in the hyperedge.  We say that the closest
vertex is where the hyperedge \emph{starts} and the others are where
it \emph{ends}.  The \emph{rooted degree} of a vertex $i$ in a planted
hyperforest is the number of hyperedges starting at $i$.  Notice that
the rooted degree of a vertex is one less than its valence (i.e. the
number of hyperedges that contain $i$) unless that vertex happens to
be a root, in which case the two numbers are equal.  Finally, to
establish notation, let $\PHF_n^k$ denote the set of all planted
hyperforests on $[n]$ with exactly $k$ components and let $\RHT_n$
($=\PHF_n^1$) denote the set of all rooted hypertrees on $[n]$.

We establish a pair of identities involving summations of symmetric
polynomials over planted hyperforests and (unrooted) hypertrees,
respectively.  We begin by defining a particular weight function for
rooted hyperforests that is similar to the block weights and partition
lattice functions consider above.

\begin{defn}[Weight of a planted hyperforest]\label{def:rooted-weight}
Let $\tau$ be a planted hyperforest.  Define the \emph{weight of an
edge} in $\tau$ as $(e-1)^{e-2}$ where $e$ is its size.  Define the
\emph{weight of the vertex $i$} as $x_i^{\deg i}$ where $\deg i$ is
its rooted degree.  The \emph{weight of $\tau$} is the product of all
of its edge and vertex weights.  Thus, \[ \wt(\tau) = \left(
\prod_{e\in \es(\tau)} (e-1)^{e-2} \right)\cdot \left( \prod_{i\in
\vs(\tau)} x_i^{\deg i}\right).\]
\end{defn}

The main combinatorial result in this article is an assertion that the
sum over all planted hyperforests of their weights is equal to a very simple
symmetric function.  This formula is the key step in establishing
Theorem~\ref{mainthm:owh}.

\begin{thm}[A hyperforest identity]\label{thm:hyperforest}
When the weight of a planted hyperforest is defined as above, the sum
of these weights over all planted hyperforests on $[n]$ with exactly
$k$ components is a function of the linear symmetric function.  More
specifically,
\begin{equation}\label{eq:rhfid}
 \sum_{\tau \in \PHF_n^k} \wt(\tau) = \binom{n-1}{k-1} h_{n-k}(s_{[n]})
\end{equation}
which specializes to 
\begin{equation}\label{eq:rhtid}
 \sum_{\tau \in \RHT_{n+1}} \wt(\tau) = h_{n}(s_{[n+1]}) = s_{[n+1]}
 (s_{[n+1]} + n)^{n-1}
\end{equation}
when summing these weights over rooted hypertrees on $[n+1]$.
\end{thm}

\begin{proof}
The proof is by induction on the ordered pairs $(n,k)$, ordered
lexicographically.  Since the base case $(1,1)$ is easily checked,
suppose that the result has been established for all ordered pairs
less than $(n,k)$.  The inductive step splits into two cases depending
on whether $k=1$ or $k>1$.  When $k>1$, then the sum on the lefthand
side can be rearranged according to the partition of $[n]$ determined
by the connected components of the hyperforest under consideration.
Once a partition $\sigma$ in $\Pi_n^k$ has been fixed, the
hyperforests associated to it are found by picking a rooted hypertree
on each block of $\sigma$.  Because these choices can be made
independent of one another, and because the weight of a planted
hyperforest is the product of the weights of its rooted hypertrees,
the sum of the hyperforest weights over the hyperforests associated
with this particular partition is a product of sums that only focus on
one block of $\sigma$ at a time.  In addition, if $\tau$ is a block of
$\sigma$ then the inductive hypothesis can be applied to the sum of
rooted hypertree weights over each all rooted hypertrees on $\tau$
since the size of $\tau$ must be strictly less than $n$.  Thus we have
\[ 
\sum_{\mu \in \PHF_n^k} \wt(\mu) = \sum_{\sigma\in \Pi_n^k}
\prod_{\tau \in \blocks{\sigma}} h_{|\tau|-1}(s_\tau)\ , 
\] 
which is equal to $\binom{n-1}{k-1} h_{n-k}(s_{[n]})$ by the $3^{rd}$
partition identity, Corollary~\ref{cor:part3}.

When $k=1$, the sum on the lefthand side is only over rooted
hypertrees and the approach is slightly different.  In order for the
algebraic manipulations to be easier to follow, we consider the case
$(n+1,1)$ rather than $(n,1)$.  If $\mu$ is a rooted hypertree on
$[n+1]$, then there is a planted hyperforest on $n$ vertices
associated with $\mu$ which can be found by removing its root and all
of the hyperedges containing it.  The resulting components are
naturally rooted by selecting the unique vertex in each that was in a
common hyperedge with the original root.  (See
Figure~\ref{fig:hypercorresp}.)

\begin{figure}[hptb]
\includegraphics{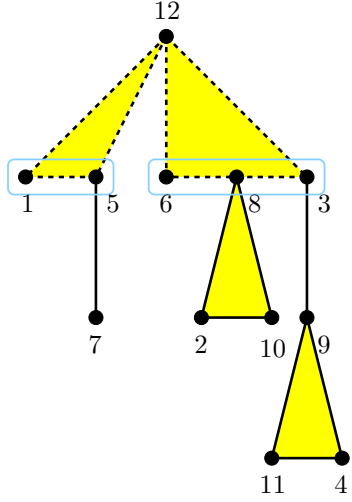}
 \caption{A rooted hypertree on $[12]$ with its associated planted
 hyperforest on $[11]$.  One can pass from the hyperforest back to the
 hypertree using the partition of $\{1, 3, 5, 6, 8\}$ indicated by the
 lightly colored boxes.}
 \label{fig:hypercorresp}
\end{figure}

As in the previous case, we proceed by rearranging a sum and applying
the inductive hypothesis.  Consider the portion of the lefthand side
of equation~\ref{eq:rhfid} that only sums over those rooted hypertrees
on $[n+1]$ that have the vertex $n+1$ as their root, the root has
rooted degree exactly $i$, and the number of connected components in
the planted forest on $[n]$ found by removing the root and all its
hyperedges is exactly $k$.  The collection of all rooted hypertrees
having all three properties can be found as follows.  First, pick an
arbitrary planted forest $\mu$ with exactly $k$ components on the
vertex set $[n]$.  Next, pick a partition $\sigma$ of the $k$ roots of
$\mu$ that has exactly $i$ blocks.  Finally, for each block in
$\sigma$, add a hyperedge containing the appropriate roots of $\mu$
together with the vertex $n+1$.  (This process is illustrated in
Figure~\ref{fig:hypercorresp} where $k=5$ and the two blocks of the
partition are indicated by the lightly colored loops.)  Because of the
independence of the choices involved, the sum of the hypertree weights
over this restricted set of rooted hypertrees can be rewritten as
follows.
\[ 
\left( \sum_{\sigma\in \Pi_k^i} \prod_{\tau \in\blocks(\sigma)}
|\tau|^{|\tau|-1} \right) \left( \sum_{\mu \in \PHF_n^k} \wt(\mu)
\right) x_{n+1}^i 
\]
The first factor is the contribution of the edge weights of the
hyperedges containing the root (there is a slight shift since the size
of the block in $\sigma$ is one less than the size of the edge
containing the root), the third factor is the contribution of the
vertex weight of the root itself, and the second factor is the
contribution of all of the other vertex and edge weights combined.  By
the $2^{nd}$ partition identity, from Corollary~\ref{cor:part2}, the
first factor is equal to $\binom{k-1}{i-1} k^{k-i}$, and by the
inductive hypothesis the second factor is equal to $\binom{n-1}{k-1}
h_{n-k}(s_{[n]})$.  Using the fact that
\[
\binom{k-1}{i-1} \binom{n-1}{k-1} =
\binom{n-1}{i-1,k-i,n-k} = \binom{n-1}{i-1} \binom{n-i}{k-i},
\] 
the sum of these weights can be rewritten as
\[ 
\binom{n-1}{i-1}
\binom{n-i}{k-i} k^{k-i} h_{n-k}(s_{[n]}) x_{n+1}^i.
\] 
A similar formula holds when any other vertex is selected as root.

Finally, it only remains to eliminate the dependence of this sum on
$k$, $i$, and the choice of a root vertex.  One form of Abel's
identity, Equation~\ref{eq:abel1}, is
\[
\sum_k \binom{n-i}{k-i}
h_{n-k}(s_{[n]}) k^{k-i} = (s_{[n]}+n)^{n-i},
\] 
where we have substituted $s_{[n]}$ for $x$, $i$ for $y$, $n-i$ for
$n$, $n-k$ for $i$, and $k-i$ for $j$.  Thus, summing over all
possible values of $k$ we find that the sum of the weights over all
hypertrees with root vertex $n+1$ and root degree $i$ is equal to
$\binom{n-1}{i-1} (s_{[n]}+n)^{n-i} x_{n+1}^i$.  Summing over all
possible values of $i$, factoring out a single $x_{n+1}$ and applying
the binomial theorem we find that the sum of the weights of all
hypertrees with root $n+1$ is $x_{n+1} (s_{[n]} + n + x_{n+1})^{n-1} =
x_{n+1} (s_{[n+1]} + n)^{n-1}$.  Finally, summing over the possible
choices of the root vertex gives $s_{[n+1]} (s_{[n+1]}+n)^{n-1}$ which
is precisely $h_n(s_{[n+1]})$ as desired.  This completes the
inductive step when $k=1$ and the proof.
\end{proof}

As is often the case, we can use the identity established for rooted
hypertrees to establish a similar identity for \emph{unrooted}
hypertrees.

\begin{defn}[Weight of a hyperforest]
If $\tau$ is an (unplanted) hyperforest, then the \emph{weight of an
edge} in $\tau$ is as before (i.e. $(e-1)^{e-2}$ where $e$ is its
size) but the \emph{weight of the vertex $i$} is now
$x_i^{\mbox{\scriptsize \val}(i) -1}$ where $\val(i)$ is its valence.
The \emph{edge weight} of $\tau$ is the product of the weights of its
edges,
\[
\ewt(\tau)= \prod_{e \in \es(\tau)} (e -1)^{(e-2)}.
\]
The \emph{vertex weight} of $\tau$ is the product of the weights of
its vertices,
\[
\vwt(\tau) = x_1^{\mbox{\scriptsize \val}(1) - 1} 
             x_2^{\mbox{\scriptsize \val}(2)-1} 
             \cdots x_n^{\mbox{\scriptsize \val}(n)-1},
\]
which is closely related to the ``degree sequence monomial'' of a
planted forest.  The \emph{weight of $\tau$} is
\[
\wt(\tau) = \ewt(\tau) \cdot \vwt(\tau), 
\]
so it is still the product of $\tau$'s edge and vertex weights.
\end{defn}

Notice that these definitions almost completely agree with those given
in Definition~\ref{def:rooted-weight}.  In particular, if $\tau$ is a
rooted hypertree on $[n]$ rooted at vertex $i$, then its weight when
viewed as a rooted hypertree is $x_i$ times its weight when viewed as
an unrooted hypertree.  This is because $\val(j) -1 = \deg j$ except
at the root.  As a consequence, if $\tau$ is a hypertree and $\tau_i$
is the hypertree $\tau$ with root vertex $i$, then the sum of the
rooted weights of the $\tau_i$ is $s_{[n]}$ times the weight of
$\tau$.

\begin{thm}[A hypertree identity]\label{thm:comb-main}
When the weight of a hypertree is defined as above, the sum of these
weights over all hypertrees on $[n+1]$ is a function of the linear
symmetric function.  More specifically,
\begin{equation}\label{eq:htid}
 \sum_{\tau \in \hyper_{n+1}} \wt(\tau) = (s_{[n+1]} + n)^{n-1} 
\end{equation}
\end{thm}

\begin{proof} 
By the remark made above, the sum of the rooted weights over all
rooted hypertrees on $[n+1]$ is equal to $s_{[n+1]}$ times the sum of
the weights over all (unrooted) hypertrees on $[n+1]$.  The formula
now follows from the second formula in the statement of
Theorem~\ref{thm:hyperforest}.
\end{proof}

Here are two small examples that illustrate the equality.

\begin{exmp}[Low dimensional cases]
When $n=2$, the poset $\hyper_3$ has four elements and the four terms
on the lefthand side are $1\cdot x_1$, $1\cdot x_2$, $1\cdot x_3$ and
$2\cdot 1$.  Thus the total is $x_1+x_2+x_3+2$ as predicted by the
righthand side.

There are four types of hypertrees in $\hyper_4$, and $29$ elements
total in the poset.  There is one element consisting of a single
hyperedge (Type A); twelve elements consisting of two hyperedges (Type
B); twelve elements that are simplicial arcs (Type C); and four
elements that are simplicial trees containing a central trivalent
vertex (Type D).  (These are drawn in Figure~3 of \cite{McMe04}.)
Starting on the lefthand side, we have that the element of type $A$
contributes $3^2$, the elements of type $B$ add $6(x_1+x_2+x_3+x_4)$,
the elements of type $C$ add
$2(x_1x_2+x_1x_3+x_1x_4+x_2x_3+x_2x_4+x_3x_4)$ and the elements of
type $D$ add $x_1^2+x_2^2+x_3^2+x_4^2$.  So for $n=3$ the total is
$(x_1+x_2+x_3+x_4+3)^2$ as expected.
\end{exmp}

\section{Euler, McCullough and Miller}\label{sec:euler}

Our Euler characteristic results follow from a computation using the
action of outer Whitehead automorphism groups on certain contractible
complexes introduced by McCullough and Miller \cite{McMi96}.  We
recall relevant facts about Euler characteristics of groups and about
these McCullough-Miller complexes.

\subsection{Euler characteristics}

Because the homotopy type of an aspherical complex is determined by
its fundamental group, it makes sense to define the Euler
characteristic of a group $G$ as the Euler characteristic of any
finite $K(G,1)$.  However, a group does not have to be of finite type
to have a well defined Euler characteristic.  A group $G$ is of
\emph{finite homological type} if $G$ has finite $\vcd$ and for every
$i$ and every $G$-module $M$ that is finitely generated as an abelian
group, $H_i(G,M)$ is finitely generated.  This holds, for example, if
$G$ is VFL, that is, has a finite index subgroup $H$ where $\Z$ admits
a finite free resolution as a trivial $\Z H$-module.  The Euler
characteristic can be defined for any group of finite homological
type.

To complete the proofs of Theorems~\ref{mainthm:owh},
\ref{mainthm:aut}, and \ref{mainthm:palindromic}, we need some basic
facts about Euler characteristics of groups.  The following results
are parts $(b)$ and $(d)$ of Theorem~$7.3$ in Chapter~IX of
\cite{Br94}.

\begin{thm}\label{thm:extensions}
Let $G$ be a group of finite homological type.
\begin{enumerate}
\item \label{finite-index}
If  $H$  is a
subgroup of finite index in $G$ then $\euler(H) = \euler(G)\cdot
[G:H]$.
\item \label{short-exact}
If $1\to K \to G \to Q \to 1$ is a short exact sequence of groups (with
$K$ and $Q$ of
finite homological type) then $\euler(G) = \euler(K) \cdot \euler(Q)$.
\end{enumerate}
\end{thm}

Our approach to computing the Euler characteristics of Whitehead
automorphism groups is via the standard formula for equivariant Euler
characteristics (see \cite{Br74} and \S IX.7 of \cite{Br94}).

\begin{thm} \label{thm:euleraction}
Let $X$ be a cocompact, contractible $G$-complex, and let $\mathcal E$
be a set of representatives of the cells of $X$ mod $G$.
Then 
\[
\euler(G) = \sum_{\sigma \in {\mathcal E}}
\left(-1\right)^{\mbox{\scriptsize{\rm dim}}(\sigma)}
\euler\left(G_\sigma\right)\ .
\] 
In particular, if $G = \ast_{i=1}^n
G_i$ is a free product of groups of finite homological type, then the
associated action of $G$ on its Bass-Serre tree shows that 
\[
\euler(G) = \euler(G_1) + \cdots + \euler(G_n) - (n - 1)\ .
\]
\end{thm}

Euler characteristics of groups have been computed in a number of
interesting cases.  For example, the Euler characteristics of many
arithmetic groups can be expressed in terms of zeta-functions, an
interesting example being 
\[ 
\euler\left(\mbox{SL}_2(\Z[1/p])\right) =
\frac{p-1}{12} 
\] 
(see the survey \cite{Se79}).  The Euler characteristics of mapping
class groups of surfaces have also been explicitly computed
(\cite{HaZa86, Pe88}).  In general Euler characteristics of families
of groups do not have tidy formulas.  For example, a generating
function for $\euler\left(\Out(F_n)\right)$ is presented in
\cite{SmVo87}, but the situation is sufficiently complicated that no
closed form solution is given.  Computer computations are used to find
$\euler(\Out(F_n))$ up to $n=12$, where
\[ 
\euler(\Out(F_{12})) = -
\frac{375393773534736899347}{2191186722816000}\ . 
\] 
For further
information on Euler characteristics of groups read Chapter IX of
\cite{Br94}.

\subsection{McCullough-Miller complexes}

Given a free product $G = G_1 \ast \cdots \ast G_n$, McCullough and
Miller construct a $\OWh(G)$ complex that can be viewed as a space of
certain actions of $G$ on simplicial trees.  At the heart of their
construction is the poset of hypertrees.

\begin{figure}[hptb]
 \includegraphics{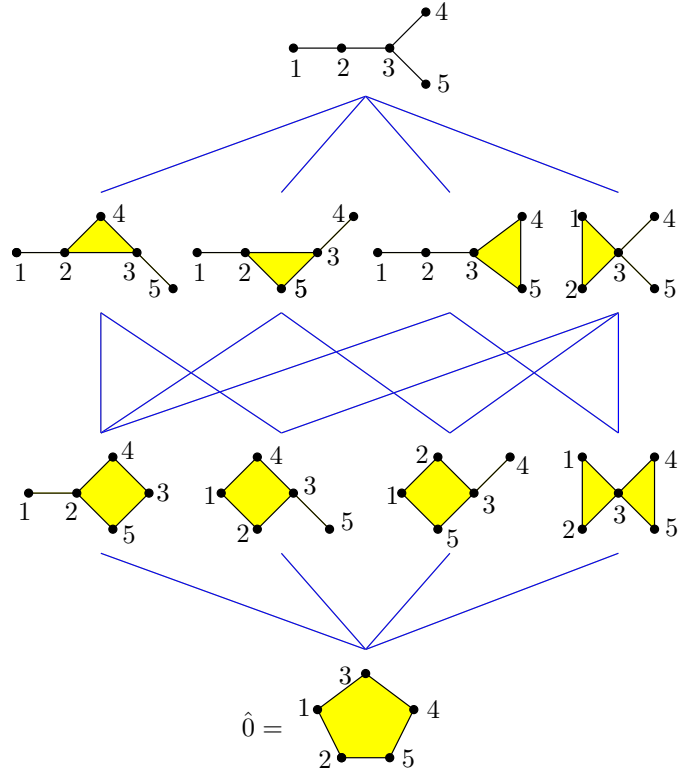}
 \caption{The Hasse diagram of the closed interval between $\hat{0}$ and 
one of the $125$ elements of maximal rank  in $\hyper_5$.}
 \label{fig:interval}
 \end{figure}

\begin{defn}[The hypertree poset]
The set of all hypertrees on $V=[n]$ forms a poset, $\hyper_n$, where
$\tau \leq \tau'$ if each edge of $\tau'$ is contained in an edge of
$\tau$.  We write $\tau < \tau'$ if $\tau \leq \tau'$ but $\tau \neq
\tau'$.  The poset $\hyper_n$ is a graded poset and the hypertrees at
height $i$ are precisely those hypertrees with $i+1$ edges.  Notice
that $\hyper_n$ contains a unique minimal element (with height $0$)
that has only a single edge containing all of $[n]$.  This is the
\emph{nuclear vertex} and we denote it by $\hat 0$.  In
Figure~\ref{fig:interval} we show the Hasse diagram for the interval
$[\hat{0}, \tau] \subset \hyper_5$ where $\tau$ is a maximal element
of height $3$.  For convenience we have used convex polygons (edges,
triangles, squares, etc.) to indicate the hyperedges.  Since
$\hyper_5$ contains $311$ elements, $125$ of which are maximal
elements, we do not show the entire poset.  If $\tau \in \hyper_n$ let
$\hyper_{n,\ge \tau}$, and $\hyper_{n,> \tau}$, denote the subposets
induced by elements which are above, or strictly above, $\tau$.
\end{defn}

An easy induction establishes the following well-known fact.

\begin{lem}\label{lem:htrank}
For any hypertree $\tau$, 
\[
\mbox{\rm Height}(\tau) = \#E(\tau) - 1 =
\sum_{v \in V(\tau)} (\val(v) - 1)\ ,
\] 
where $\val(v)$ is the valence
of $v$, i.e. the number of hyperedges containing $v$.
\end{lem}

\noindent 
We refer to the value multiply described in Lemma~\ref{lem:htrank} as
the \emph{rank} of the hypertree $\tau$, and denote it $\rk(\tau)$.

The main result of \cite{McMi96} can  be formulated as follows.

\begin{thm} \label{thm:funddomain}
Given $G = G_1 \ast \cdots \ast G_n$ there is a contractible
$\OWh(G)$-complex $\MM_n$ with $\OWh(G)\backslash \MM_{n} \simeq
|\hyper_n|$.
\end{thm}

\noindent Here $|\hyper_n|$ denotes the geometric realization of
$\hyper_n$.  The reader wishing to understand the connection between
this description and the original description as given in
\cite{McMi96} can consult \cite{McMe04}.  Geometrically, the complex
$\MM_n$ contains an imbedded copy of $|\hyper_n|$ that projects
isomorphically to the quotient under the action, and McCullough and
Miller give a clean description of simplex stabilizers.

\begin{thm} \label{thm:stabilizers}
The stabilizer of $\tau \in |\hyper_n|$ under the action of $\OWh(G)$
is
\[ 
\mbox{\rm Stab}(\tau) = \bigoplus_{i=1}^n 
\left(G_i\right)^{\mbox{\scriptsize \val}(i)-1}\ . 
\] 
Further, if $\tau < \tau' < \cdots < \tau''$ is the chain
corresponding to a simplex $\sigma \in |\hyper_n|$, then
$\mbox{Stab}(\sigma) = \mbox{Stab}(\tau)$.
\end{thm}

If the $G_i$ are all of finite homological type, then the simplex
stabilizers in $\MM_n$ under the action of $\OWh(G)$ are of finite
homological type.  Using the fact that $\MM_n$ is a finite
dimensional, contractible, cocompact $\OWh(G)$-complex, we get:

\begin{cor}\label{cor:finhomtype}
If $G = \ast_{i=1}^n G_i$, and each of the $G_i$ is of
finite homological type, then $\OWh(G)$ is of finite homological
type.
\end{cor}

\begin{rem}
Even if $G$ is of finite homological type, $\Aut(G)$ may not be.  For
example, the Baumslag-Solitar group $\mbox{BS}(2,4)$ is of finite
type, yet $\Aut(\mbox{BS}(2,4))$ is not finitely generated
\cite{CoLe83}.  Corollary~\ref{cor:finhomtype} establishes that
$\euler(\OWh(G))$ and $\euler(\Wh(G))$ exist when $G$ is a free
product of groups of finite homological type, but the example of
$\mbox{BS}(2,4)$ illustrates that in general $\euler(\Aut(G))$ may not
exist.  Thus one cannot hope to establish a result such as
Theorem~\ref{mainthm:aut} without some reasonable hypotheses such as
``all $G_i$ finite.''
\end{rem}

We record a few useful facts about the combinatorics of
$\hyper_n$ from \cite{McMe04}.

\begin{lem} \label{lem:intervals}
Let $\tau$ be a hypertree on the set $[n]$.
 The subposet $\hyper_{n,\ge \tau}$ decomposes as a direct product of
  hypertree posets:
\[
\hyper_{n,\ge \tau} = \prod_{j \in \es(\tau)} \hyper_{j}\ .
\]
\end{lem}

Recall that the \emph{reduced} Euler characteristic of a
complex $K$ is $\widetilde{\euler}(K) = \euler(K) - 1$, that
is, it is the Euler characteristic based on reduced homology.

\begin{thm}\label{thm:mobius}
For every $n$, 
\[
\widetilde{\euler}(|\hyper_{n+1, > \hat{0}}|)
= (-1)^n n^{n-1}\ .
\]
\end{thm}

Theorem~\ref{thm:mobius} and Lemma~\ref{lem:intervals} combine to give
the following formula.  The only tricky matter is establishing the
sign.  See \cite[Lemma~4.6]{McMe04} for the details.

\begin{cor} \label{cor:asclinks}
The reduced Euler characteristic of $|\hyper_{n,>\tau}|$
is 
\begin{eqnarray*}
\widetilde{\euler}(|\hyper_{n,>\tau}|) 
        &=& (-1)^{n + |\mbox{\scriptsize \edges}(\tau)|} 
            \prod_{j \in \es(\tau)} (j-1)^{(j-2)}\\ 
        &=& (-1)^{n + |\edges(\tau)|} \cdot \ewt(\tau)\ .\\ 
\end{eqnarray*}
\end{cor}

\section{Proofs and Computations}\label{sec:main}
In this final section we complete the proofs of
Theorems~\ref{mainthm:owh} and \ref{mainthm:aut}.  First, we express
the Euler characteristic of $\OWh(G)$ in terms of a sum over
hypertrees.

\begin{lem}\label{lem:euler-owh}
Let $G = \ast_{i=1}^n G_i$, where each $G_i$ is of
finite homological type.  Then
\[
\euler(\OWh(G)) = (-1)^n \sum_{\tau\in \hyper_n}  
  (-1)^{\mbox{\scriptsize \rk}(\tau)} \cdot 
\euler(\mbox{\rm Stab}(\tau))\cdot \ewt(\tau)\ .
\]
\end{lem}

\begin{proof}
By Theorem~\ref{thm:euleraction} we know 
\[
\euler(\OWh(G)) = \sum_{\sigma \in |\hyper_n|} 
    (-1)^{\mbox{\scriptsize{dim}}(\sigma)}
\euler(\mbox{\rm Stab}(\sigma))\ .
\]
If $\sigma = \tau < \tau' < \cdots < \tau''$, 
 the stabilizer of $\sigma$ is the same as the stabilizer
of $\tau$.  Thus we may reorganize this summation as
\begin{eqnarray*}
\euler(\OWh(G)) &=& \sum_{\tau \in \hyper_n} 
 \left\{ \euler(\mbox{\rm Stab}(\tau)) \sum_{\sigma \in |\hyper_{n, \ge \tau}|}
 (-1)^{\mbox{\scriptsize{dim}}(\sigma)}\right\}\\ 
&=& \sum_{\tau \in \hyper_n} 
 \left\{ \euler(\mbox{\rm Stab}(\tau)) \left(1- \sum_{\sigma \in |\hyper_{n, > \tau}|}
 (-1)^{\mbox{\scriptsize{dim}}(\sigma)}\right)\right\}\\ 
&=& \sum_{\tau \in \hyper_n} 
 \euler(\mbox{\rm Stab}(\tau)) \left(
 -\widetilde{\euler}(|\hyper_{n, >\tau}|\right)\ .\\
\end{eqnarray*}
But up to sign, the reduced Euler characteristic of $|\hyper_{n,> \tau}|$ is just
$\ewt(\tau)$ (Corollary~\ref{cor:asclinks}).  Inserting
the formula from Corollary~\ref{cor:asclinks} we get
\begin{eqnarray*}
\euler(\OWh(G)) 
&=&  \sum_{\tau \in \hyper_n} 
 \euler(\mbox{Stab}(\tau)) \left((-1)^{n + 1 + |\edges(\tau)|}
  \cdot \ewt(\tau)\right)\\ 
&=& (-1)^n \sum_{\tau \in \hyper_n} 
 \euler(\mbox{Stab}(\tau)) \cdot (-1)^{\mbox{\scriptsize \rk}(\tau)} \cdot 
  \ewt(\tau)\ ,
\end{eqnarray*}
which is the formula claimed in the statement of the lemma.
\end{proof}

Theorem~\ref{mainthm:owh} is now essentially established.  
Exploiting Corollary~\ref{cor:finhomtype} we can give it a slightly
stronger phrasing.

\setcounter{mainthm}{0}
\begin{mainthm}
If $G = G_1 \ast \cdots \ast G_n$ is a free product of groups of finite
homological type, then the Euler characteristic of the group of outer
Whitehead automorphisms is $\euler(\OWh(G))=\euler(G)^{n-2}$.  As a
consequence, the Euler characteristic of the group of Whitehead
automorphisms is $\euler(\Wh(G))=\euler(G)^{n-1}$.
\end{mainthm}

\begin{proof}
Because $\mbox{Stab}(\tau) =\bigoplus_{i=1}^n
\left(G_i\right)^{\mbox{\scriptsize \val}(i)-1}$ by
Theorem~\ref{thm:stabilizers}, and using the equivalence of the
multiple definitions of rank given in Lemma~\ref{lem:htrank}, we see
that
\[
 \euler(\mbox{Stab}(\tau)) \cdot (-1)^{\mbox{\scriptsize \rk}(\tau)}
 = \prod ( - \euler(G_i))^{\mbox{\scriptsize \val}(i) - 1}\ .
\]
So we can rewrite the formula from Lemma~\ref{lem:euler-owh} as
\[
\euler(\OWh(G)) 
  =   (-1)^n \sum_{\tau \in \hyper_n} 
 \left[ \prod (- \euler(G_i))^{\mbox{\scriptsize \val}(i) - 1} \right]
 \cdot  \ewt(\tau)\ .
\]
Thus---aside from the $(-1)^n$ term---the Euler characteristic of
$\OWh(G)$ is given by the left side of hypertree
identity~\ref{eq:htid}, with ``$-\euler(G_i)$'' being substituted for
the variable ``$x_i$''.  Making the same substitution on the right
side of the hypertree identity and reinserting the $(-1)^n$ term shows
that
\begin{eqnarray*}
\euler(\OWh(G))&=& (-1)^n \left(-\euler(G_1) - \euler(G_2) - \cdots - \euler(G_n)
+ (n-1)\right)^{n-2}\\ 
&=& (-1)^n \left(-\euler(G)\right)^{n-2} = \euler(G)^{n-2}.
\end{eqnarray*}
Thus we have the formula for the group of outer Whitehead automorphisms.

Theorem~\ref{thm:extensions}.\ref{short-exact} applied to the short
exact sequence 
\[
1 \to G \to \Wh(G) \to \OWh(G) \to 1
\]
 shows that
$\euler(\Wh(G)) = \euler(G)\cdot \euler(\OWh(G))= \euler(G)^{n-1}$.
\end{proof}

\begin{mainthm}
If $G = G_1 \ast \cdots \ast G_n$ is a free product of finite groups and
$\Omega \subset \Sigma_n$ is the subgroup of the symmetric group given
by permuting isomorphic factors $G_i$, then the Euler characteristics
of the Fouxe-Rabinovitch automorphism group of $G$, the full
automorphism group of $G$ and the outer automorphism group of $G$ are
as follows.
\[\begin{array}{rcl}
\euler(\FR(G)) &=& \euler(G)^{n-1}  \displaystyle\prod_{i=1}^n |\Inn(G_i)|\\
\euler(\Aut(G)) &=& \euler(G)^{n-1} |\Omega|^{-1} 
    \displaystyle\prod_{i=1}^n |\Out(G_i)|^{-1}\\
\euler(\Out(G)) &=& \euler(G)^{n-2} |\Omega|^{-1} 
    \displaystyle\prod_{i=1}^n |\Out(G_i)|^{-1}\\
\end{array}\]
\end{mainthm}

\begin{proof}
The identities listed above are derived from Theorem~\ref{mainthm:owh}
using finite-index subgroups  and
short exact sequences.
To calculate $\euler(\FR(G))$ we use the fact that $\FR(G)$ is a
finite-index subgroup of $\Wh(G)$.  Recall that if each $G_i$ is finite,
then $\Aut(G) = \Sigma \Aut(G)$ and so $\Wh(G)$ is the kernel
of the map $f:\Aut(G)\to \Out(\prod_i G_i)$ and $\FR(G)$ is the
kernel of the map $g:\Aut(G) \to \Aut(\prod_i G_i)$.  Moreover, the
map $f= h \circ g$ where $h$ is the map $\Aut(\prod_i G_i)\to
\Out(\prod_i G_i)$.  Since the kernel of $h$ is $\prod_i \Inn(G_i)$,
the index of $\FR(G)$ in $\Wh(G)$ is $\prod_i |\Inn(G_i)|$.  This index
is finite when each $G_i$ is finite, and the first identity follows
by Theorem~\ref{thm:extensions}.\ref{finite-index}.

To calculate $\euler(\Aut(G))$ we use the decomposition of $\Aut(G)$
noted by McCullough and Miller \cite[p.~48]{McMi96}:
\[
\Aut(G) = FR(G) \rtimes \left(\left(\prod_i \Aut(G_i)\right) \rtimes
\Omega\right)\ .
\]
From this decomposition it is clear that the size of $g(\Aut(G))$---and 
the index of $\FR(G)$ in $\Aut(G)$---is $|\Omega| \prod_i
|\Aut(G_i)|$.  Dividing by the index of $\FR(G)$ in $\Wh(G)$ we see
that the index of $\Wh(G)$ in $\Aut(G)$ is $|\Omega| \prod_i
|\Out(G_i)|$.  This index is finite and the second identity follows,
again by \ref{thm:extensions}.\ref{finite-index}.  Finally, from the
short exact sequence $1 \to G \to \Aut(G) \to \Out(G) \to 1$ we find
that $\euler(\Out(G)) = \euler(G)^{-1} \euler(\Aut(G))$ by
Theorem~\ref{thm:extensions}.\ref{short-exact} which completes the
proof.
\end{proof}

\begin{exmp}[The factors are all finite cyclic]
When $G$ is a free product of finite cyclic groups, the values listed
in Theorem~\ref{mainthm:aut} are particularly easy to compute.
If $G_i$ is finite cyclic with order $g_i$, then $\euler(G_i) =
(g_i)^{-1}$ and $\euler(G) = (\sum_i (g_i)^{-1} -n +1)$.  Moreover,
$G_i$ abelian implies $|\Inn(G_i)| = 1$ and $G_i$ cyclic means that
$|\Aut(G_i)| = |\Out(G_i)| = \phi(g_i)$ where $\phi$ is Euler's
totient function.  Recall that $\phi$ can be defined be setting
$\phi(p^k)=p^k-p^{k-1}$ for all prime powers and then extending it by
$\phi(m\cdot n) = \phi(m)\cdot \phi(n)$ wherever $m$ and $n$ are
relatively prime.  Finally, $\Omega$ is simply a product of symmetric
groups whose indices count the number of times each particular value
occurs among the $g_i$.  As a concrete illustration, when $G= \Z_2
\ast \Z_2 \ast \Z_2 \ast \Z_2 \ast \Z_3 \ast \Z_3 \ast \Z_3$, then
$\euler(G) = (\frac42 + \frac33 -7+1) = -3$, $|\Omega| = (4! \cdot
3!)$, $\phi(2)=1$ and $\phi(3)=2$.  Thus $\euler(\Aut(G)) =
\frac{(-3)^6}{4! 3!} (\frac11)^4 (\frac12)^3 = \frac{81}{128}$.
\end{exmp}

\begin{exmp}[The factors are all infinite cyclic]
If each $G_i \simeq \Z$, then $\Wh(\ast_{i=1}^n G_i)$
is the subgroup of $\Aut(F_n)$ 
called the \emph{pure symmetric automorphism group},
commonly denoted $P\Sigma_n$.  
The Euler characteristic of this group was computed in
\cite{McMe04}, but the answer also follows from Theorem~\ref{mainthm:owh}.
Using the fact that $\euler(F_n) = 1-n$ one gets
\[
\euler(\Wh(\underbrace{\Z \ast \cdots \ast \Z}_{n\mbox{\scriptsize 
                \ free factors}})) 
        = \euler(P\Sigma_n) = (1-n)^{(n-1)}\ .
\]
However, unlike the case where the $G_i$ are finite, $\Wh(F_n)$ is of
infinite index in $\Aut(F_n)$.  In fact, the cohomological dimension
of $\Wh(F_n)$ is $n-1$ \cite{Co89}, precisely half the $\vcd$ of
$\Aut(F_n)$ \cite{CuVo86}.
\end{exmp}

\section{Palindromic automorphisms}\label{sec:palindromic}

The palindromic automorphism group $\Pi A(F_n)$ of a free group
$F_n = F(a_1 , \ldots,  a_n )$
was introduced by Collins \cite{Co95}.  A \emph{palindrome}
is a word $w_1 w_2 \ldots w_n$ in $F_n$ that is the same as its
\emph{reverse} $w_n \ldots w_2 w_1$, and the group $\Pi A(F_n)$ consists
of all automorphisms of $F_n$ that send each generator of $F_n$ to
a palindrome.  Collins derived a presentation for $\Pi A(F_n)$
using the following generators:
\begin{itemize}
\item {\em Palindromic Whitehead moves} $(a_i||a_j)$, $i \not = j$,
which send $a_i \mapsto a_j a_i a_j$
and fix all other generators $a_k$.
\item {\em Factor automorphisms} $\sigma_{a_i}$ which send
$a_i \mapsto a_i^{-1}$ and fix all other generators $a_k$.
\item {\em Permutation automorphisms} corresponding to elements of
the symmetric group $\Sigma_n$ which permute the $a_1, \ldots, a_n$
among themselves.
\end{itemize}
The subgroup of $\Pi A(F_n)$ generated just by the palindromic
Whitehead moves $(a_i||a_j)$ is the {\em elementary palindromic
automorphism group} $E \Pi A(F_n)$ and the subgroup generated by the
palindromic Whitehead moves and the factor automorphisms is the {\em
pure palindromic automorphism group} $P \Pi A(F_n)$.

In \cite{GlJe00} Glover and Jensen observe that $\Pi A(F_n)$ is a
centralizer of an involution: $\Pi A(F_n) = C_{\Aut(F_n)}(\sigma_n)$
where $\sigma_n$ is the automorphism that sends each $a_i$ to
$a_i^{-1}$.  This allows them to show that $E \Pi A(F_n)$ is torsion
free, $\mbox{cd}\left(E\Pi A(F_n)\right) = n-1$, as well as to find
information on the cohomology of $\Pi A(F_n)$.

\vspace{10pt}

In this section we calculate the Euler characteristic of $\Pi A(F_n)$.
We thank Andy Miller for first observing the following useful
description of $\Pi A(F_n)$.  Let $G = \ast_{i=0}^n G_i$, where each
$G_i = \langle x_i \rangle \cong \Z_2$. Denote by $\Aut(\ast_{i=0}^n
\Z_2,x_0)$ the subgroup of $\Aut(\ast_{i=0}^n \Z_2)$ consisting of
automorphisms $\phi \in \Aut(\ast_{i=0}^n \Z_2)$ that fix the first
factor, $\phi(x_0)=x_0$.

\begin{lem}\label{lem:subgroup} There is
an isomorphism $\iota: \Aut(\ast_{i=0}^n \Z_2,x_0) \to \Pi A(F_n)$
satisfying $\iota(\FR(\ast_{i=0}^n \Z_2,x_0))= P \Pi A(F_n).$
\end{lem}

\begin{proof}
For each $i=1,\ldots, n$, let $a_i=x_0x_i$ and let $F_n$ be the
subgroup of $G$ generated by $a_1,\ldots,a_n$.  Because $F_n$ is the
subgroup of $G$ consisting of words of even length (in the $x_i$), it
is a normal subgroup of $G$ and we have a short exact sequence
\[
1 \to F_n \to G \to \Z_2 \to 1
\]
with a splitting given by sending $\Z_2$ to $G_0$.  In fact, the
subgroup of words of even length is necessarily a characteristic
subgroup of $G$ and it is a free group on the $n$ generators $a_1,
\ldots, a_n$.

Since $F_n$ is a characteristic subgroup, there is a map
\[
\iota: \mbox{Aut}(G,x_0) \to \mbox{Aut}(F_n).
\]
Let $\phi \in Ker(\iota)$.  Then $\phi(a_i)=a_i$ for all
$i=1,\ldots,n$ and $\phi(x_0)=x_0$.  Hence for each index $i$,
$\phi(x_i)=\phi(x_0x_0x_i)=\phi(x_0)\phi(x_0x_i)=
x_0\phi(a_i)=x_0a_i=x_0x_0x_i=x_i.$ Thus $\phi=1$ and we have shown
that $\iota$ is injective.

Let $\sigma_n \in \Aut(F_n)$ be the automorphism that sends each $a_i$
to $a_i^{-1}$.  We show that for each $\phi \in \Aut(G,x_0)$,
$\iota(\phi^{-1}) \sigma_n \iota(\phi) = \sigma_n$, so that
$\iota(\phi) \in \Pi A(F_n) = C_{\Aut(F_n)}(\sigma_n).$ By the
Kuro\v{s} Subgroup Theorem, there exists a sequence $x_{i_1}, \ldots,
x_{i_m}$ such that $\phi(x_i)=x_{i_1} \cdots x_{i_{m-1}} x_{i_m}
x_{i_{m-1}} \cdots x_{i_1}.$ Recall that we have defined $a_i = x_0
x_i$, and abusing notation, we set $a_0 = x_0 x_0 = 1 \in F_n$.  Note
that
\begin{eqnarray*}
[\iota(\phi^{-1}) \sigma_n \iota(\phi)](a_i) & = &
[\iota(\phi^{-1}) \sigma_n \iota(\phi)](x_0x_i) \\
& = & [\iota(\phi^{-1}) \sigma_n] (x_0\phi(x_i)) \\
& = &
[\iota(\phi^{-1}) \sigma_n]
(x_0 x_{i_1} \cdots x_{i_{m-1}} x_{i_m} x_{i_{m-1}} \cdots x_{i_1})\\
& = & [\iota(\phi^{-1}) \sigma_n]
(a_{i_1} a_{i_2}^{-1} a_{i_3} \cdots a_{i_{m-1}}^{\pm 1} a_{i_m}^{\mp 1}
a_{i_{m-1}}^{\pm 1} \cdots a_{i_1})\\
& = & \iota(\phi^{-1})
(a_{i_1}^{-1} a_{i_2} a_{i_3}^{-1} \cdots a_{i_{m-1}}^{\mp 1} a_{i_m}^{\pm
1}
a_{i_{m-1}}^{\mp 1} \cdots a_{i_1}^{-1})\\
& = & \iota(\phi^{-1})
(x_{i_1} \cdots x_{i_{m-1}} x_{i_m} x_{i_{m-1}} \cdots x_{i_1} x_0)\\
& = & x_i x_0 \\
& = & a_i^{-1}
\end{eqnarray*}
So the injection $\iota$ restricts to 
\[
\iota: \Aut(G,x_0) \to \Pi A(F_n).
\]

We show that $\iota$ is a surjection by showing that the generators of
$\Pi A(F_n)$ are in the image of $\iota$.  For each pair $i,j$ of
distinct integers from $1,\ldots, n$, let $\phi_{i,j} \in \Aut(G,x_0)$
be the automorphism sending $x_k$ to $x_0 x_k x_0$ for $k \not = i,0$
and $x_i$ to $x_jx_0x_ix_0x_j$.  Then $\iota(\phi_{i,j})=(a_i ||
a_j)$.  The permutation automorphisms of $\Pi A(F_n)$ can be realized
as maps fixing $x_0$ and permuting the $x_0 x_i x_0$.  Finally, the
factor automorphism that inverts $x_0 x_i$ and fixes all other
generators corresponds to the map where $x_i \mapsto x_0 x_i x_0$ and
all other generators are fixed.  So $\iota$ is an isomorphism.

The above calculations show that $\iota$ sends $\FR(\ast_{i=0}^k
G_i,x_0)$ to the pure palindromic automorphism group $P \Pi A(F_n)$
generated by palindromic Whitehead moves and factor automorphisms.
\end{proof}

\begin{mainthm}
If $F_n$ is a free group of rank $n$, then the Euler characteristics
of the elementary palindromic automorphism group of $F_n$, the pure
palindromic automorphism group of $F_n$, and the palindromic
automorphism group of $F_n$ are as follows.
\begin{eqnarray*}
\euler(E\Pi A(F_n)) & = & (1-n)^{n-1}\\
\euler(P\Pi A(F_n)) & = & \frac{(1-n)^{n-1}}{2^n} \\
\euler(\Pi A(F_n))  & = & \frac{(1-n)^{n-1}}{2^n \cdot n!}
\end{eqnarray*}
\end{mainthm}

\begin{proof}
Let $G$ denote the free product of $n+1$ copies of $\Z_2$ and suppose
the first $\Z_2$ is generated by $x_0$.  Conjugation by $x_0$ is the
only inner automorphism in $P\Pi A(F_n) = \FR(G,x_0) = \Wh(G,x_0)$ and
we have a short exact sequence
\[
1 \to \Z_2 \to P\Pi A(F_n) \to \OWh(G,x_0) \to 1\ .
\]
Clearly $\OWh(G,x_0) < \OWh(G)$.  However, we can conjugate each
generator of $\Wh(G)$ to form a generating set for $\OWh(G)$, each
element of which fixes $x_0$.  Hence $\OWh(G,x_0) = \OWh(G)$, and
therefore the Euler characteristic of $P\Pi A(F_n)$ is
\[
\euler(P\Pi A(F_n)) = \euler(\Z_2) \cdot \euler(\OWh(G)) = 
(1/2) \cdot ((1-n)/2)^{n-1} = \frac{(1-n)^{n-1}}{2^n}\ .
\]
The remaining assertions follow from the short exact sequences
\[
1 \to E\Pi A(F_n) \to P\Pi A(F_n) \to (\Z_2)^n \to 1
\]
and
\[
1 \to P\Pi A(F_n) \to \Pi A(F_n) \to \Sigma_n \to 1
\]
and the fact that $\euler((\Z_2)^n) = \left(\frac{1}{2}\right)^n$ and
$\euler(\Sigma_n) = \frac{1}{n!}$.
\end{proof}

\begin{rem}
The groups $\psig_n$ and $E \Pi A(F_n)$ are both torsion free groups
of cohomological dimension $n-1$, and they have very similar
presentations.  By Theorem~C they also have the same Euler
characteristic.  This might lead an optimist to think these groups are
isomorphic, or perhaps a more meek conjecture would be that they have
isomorphic homology.  However, as was noted in \cite{GlJe00},
$\psig_n$ abelianizes to a free abelian group of rank $n(n-1)$ while
$E \Pi A(F_n)$ abelianizes to an elementary abelian 2-group of rank
$n(n-1)$.
\end{rem}

Borrowing an argument of Smillie and Vogtmann---used to establish
Proposition 6.4 in \cite{SmVo87}---we can show that the kernel of
the linearization map is not of finite homological type.  As was first
pointed out in \cite{Co95}, the image of $P \Pi A(F_n)$ in the general
linear group is the congruence subgroup $\widetilde{\Gamma}_2(\Z)$
defined by the short exact sequence
\[
1 \to \widetilde{\Gamma}_2(\Z) \to GL_n(\Z) \to GL_n(\Z_2) \to 1\ .
\]
Thus if the kernel $K$ of $P\Pi A(F_n) \twoheadrightarrow
\widetilde{\Gamma}_2(\Z)$ were of finite homological type, then by
Theorem~\ref{thm:extensions}.\ref{short-exact} we would have
\[
\euler(P\Pi A(F_n)) = \euler(K) \cdot \euler(\widetilde{\Gamma}_2(\Z))
= 0
\]
since $\euler(\tilde \Gamma_2(\Z))=0$ by Harder \cite{Ha71}.  But this
contradicts Theorem~\ref{mainthm:palindromic}, so the kernel cannot
have finite homological type.  

\begin{cor}\label{cor:palindromic}
For $n \geq 3$, the kernel of the map from $P \Pi A(F_n)$ to
$GL_n(\Z)$ does not have finitely generated integral homology.
\end{cor}

As in the case of the famous kernel $\mbox{IA}_n = \ker\left[\Aut(F_n)
\twoheadrightarrow GL_n(\Z)\right]$, it would be interesting to know
what finiteness properties these kernels do enjoy.



\begin{thebibliography}{10}

\bibitem[Br74]{Br74}
Kenneth~S. Brown.
\newblock Euler characteristics of discrete groups and {$G$}-spaces.
\newblock {\em Invent. Math.}, 27:229--264, 1974.

\bibitem[Br94]{Br94}
Kenneth~S. Brown.
\newblock {\em Cohomology of groups}, volume~87 of {\em Graduate Texts in
  Mathematics}.
\newblock Springer-Verlag, New York, 1994.
\newblock Corrected reprint of the 1982 original.

\bibitem[Co89]{Co89}
Donald~J. Collins.
\newblock Cohomological dimension and symmetric automorphisms of a free group.
\newblock {\em Comment. Math. Helv.}, 64(1):44--61, 1989.

\bibitem[Co95]{Co95}
Donald~J. Collins.
\newblock Palindromic automorphisms of free groups.
\newblock In {\em Combinatorial and geometric group theory (Edinburgh, 1993)},
  volume 204 of {\em London Math. Soc. Lecture Note Ser.}, pages 63--72.
  Cambridge Univ. Press, Cambridge, 1995.

\bibitem[CL83]{CoLe83}
Donald~J. Collins and Frank Levin.
\newblock Automorphisms and {H}opficity of certain {B}aumslag-{S}olitar groups.
\newblock {\em Arch. Math. (Basel)}, 40(5):385--400, 1983.

\bibitem[CV86]{CuVo86}
Marc Culler and Karen Vogtmann.
\newblock Moduli of graphs and automorphisms of free groups.
\newblock {\em Invent. Math.}, 84(1):91--119, 1986.

\bibitem[GJ00]{GlJe00}
Henry~H. Glover and Craig~A. Jensen.
\newblock Geometry for palindromic automorphism groups of free groups.
\newblock {\em Comment. Math. Helv.}, 75(4):644--667, 2000.

\bibitem[GKP94]{Concrete}
Ronald~L. Graham, Donald~E. Knuth, and Oren Patashnik.
\newblock {\em Concrete mathematics}.
\newblock Addison-Wesley Publishing Company, Reading, MA, second edition, 1994.
\newblock A foundation for computer science.

\bibitem[Ha71]{Ha71}
G.~Harder.
\newblock A {G}auss-{B}onnet formula for discrete arithmetically defined
  groups.
\newblock {\em Ann. Sci. \'Ecole Norm. Sup. (4)}, 4:409--455, 1971.

\bibitem[HZ86]{HaZa86}
J.~Harer and D.~Zagier.
\newblock The {E}uler characteristic of the moduli space of curves.
\newblock {\em Invent. Math.}, 85(3):457--485, 1986.

\bibitem[KV93]{KrVo93}
Sava Krsti{\'c} and Karen Vogtmann.
\newblock Equivariant outer space and automorphisms of free-by-finite groups.
\newblock {\em Comment. Math. Helv.}, 68(2):216--262, 1993.

\bibitem[MM04]{McMe04}
Jon McCammond and John Meier.
\newblock The hypertree poset and the {$l\sp 2$}-{B}etti numbers of the motion
  group of the trivial link.
\newblock {\em Math. Ann.}, 328(4):633--652, 2004.

\bibitem[MM96]{McMi96}
Darryl McCullough and Andy Miller.
\newblock Symmetric automorphisms of free products.
\newblock {\em Mem. Amer. Math. Soc.}, 122(582):viii+97, 1996.

\bibitem[Pe88]{Pe88}
R.~C. Penner.
\newblock Perturbative series and the moduli space of {R}iemann surfaces.
\newblock {\em J. Differential Geom.}, 27(1):35--53, 1988.

\bibitem[Ro84]{Ro84}
Steven Roman.
\newblock {\em The umbral calculus}, volume 111 of {\em Pure and Applied
  Mathematics}.
\newblock Academic Press Inc. [Harcourt Brace Jovanovich Publishers], New York,
  1984.

\bibitem[Se79]{Se79}
J.-P. Serre.
\newblock Arithmetic groups.
\newblock In {\em Homological group theory (Proc. Sympos., Durham, 1977)},
  volume~36 of {\em London Math. Soc. Lecture Note Ser.}, pages 105--136.
  Cambridge Univ. Press, Cambridge, 1979.

\bibitem[SV87]{SmVo87}
John Smillie and Karen Vogtmann.
\newblock A generating function for the {E}uler characteristic of {${\rm
  Out}(F\sb n)$}.
\newblock {\em J. Pure Appl. Algebra}, 44(1-3):329--348, 1987.

\bibitem[St99]{StV2}
Richard~P. Stanley.
\newblock {\em Enumerative combinatorics. {V}ol. 2}, volume~62 of {\em
  Cambridge Studies in Advanced Mathematics}.
\newblock Cambridge University Press, Cambridge, 1999.
\newblock With a foreword by Gian-Carlo Rota and appendix 1 by Sergey Fomin.

\end{thebibliography}
\end{document}